\documentstyle[12pt]{article}

\newtheorem{th}{Th\'eor\`eme}

\newtheorem{prp}{Proposition}

\newtheorem{lem}{Lemme}
\newcommand{\rar}{\rightarrow}
\begin{document}
\parbox{4in}{\small  G\'eom\'etrie Alg\'ebrique / Th\'eorie de Nombres\\
              Algebraic Geometry / Number theory}

\begin{center}\large \bf Uniformit\'e des points rationnels sur tous les corps
quadratiques\\
{\normalsize \bf Dan Abramovich} \end{center}

{\small
 {\bf R\'esum\'e.} Nous pr\'ecisons  un r\'esultat de L. Caporaso, J.
Harris et B. Mazur: la conjecture de Lang entra\^{\i}ne que le nombre des
points
rationnels d'une
courbe de genre $g>1$ est born\'e uniform\'ement. Nous montrons  que cette
conjecture  entra\^{\i}ne que ce nombre est born\'e uniform\'ement pour
tous les corps quadratiques.

 \begin{center}\normalsize \bf Uniformity of rational points over all quadratic
fields
\end{center}

{\bf Abstract.} We refine a result of L. Caporaso, J. Harris and B. Mazur,
who showed that Lang's conjecture implies that the number of rational points on
a curve of genus $g>1$ is
uniformly bounded. We show that the conjecture implies that this number
is bounded uniformly  over all quadratic fields.}

Soit $K$ un corps de nombres. Selon une conjecture bien connue de S. Lang
(voir \cite{langbul}, conjecture 5.7) l'ensemble des points
$K$-rationnels d'une vari\'et\'e de type
g\'en\'eral n'est pas dense (pour la topologie de Zariski). Un travail r\'ecent
de L. Caporaso, J. Harris et B. Mazur \cite{chm}
montre que cette conjecture entra\^{\i}ne l'existence d'une borne sup\'erieure
universelle pour le
nombre des points $K$-rationnels d'une courbe lisse de genre fix\'e $g>1$.
L'ingr\'edient g\'eom\'etrique principal dans leur travail est le suivant: si
$X\rar B$ est une famille de courbes de genre $>1$, il existe $n$ tel que la
puissance relative $X^n_B$  domine une vari\'et\'e de type g\'en\'eral.

Nous allons d\'emontrer un r\'esultat plus pr\'ecis:
\begin{th}
Supposons que la conjecture de Lang soit vraie. Soit $K$ un corps des nombres
et
$g>1$ un entier. Il existe un nombre $N(K,g)$ tel que si $L$ est une
extension de degr\'e $\leq 3$ de $K$  et $C$ est une courbe lisse projective
conn\`exe, de
genre $g$ d\'efinie sur $L$ on a $$\# C(L) < N(K,g). $$
\end{th}

{\bf Remarques.}  Le cas $d=1$ est trait\'e dans \cite{chm}. Nous ne donnerons
la d\'emonstration en d\'etail que lorsque $d=2$.
Voir la remarque au d\'ebut de cette note pour le cas $d=3$.

On peut se demander s'il y a un \'enonc\'e analogue avec ``degr\'e
$\leq 3$'' remplac\'e par  ``degr\'e $\leq d$'', o\`u $d$ est un entier
quelconque. La question \`a \'et\'e pos\'ee par F. Hajir.

Je tiens a remercier E. Izadi et J.-P. Serre de m'avoir aid\'e a am\'eliorer
la qualit\'e de l'exposition et de la traduction fran\c{c}aise de ce manuscrit.

{\bf Pr\'eliminaires sur les familles de courbes.} Soit $\pi:X\rar B$ une
$K$-famille
de courbes, c'est \`a dire un morphisme projectif lisse de vari\'et\'es sur
$K$, dont les fibres sont
des courbes g\'eom\'etriquement irr\'eductibles.  On \'ecrit: $Y_n =
Sym^2(X^n_B)$; soit $B'\subset Sym^2B$,
posons $Y_{n,B'} = (Sym^2(X^n_B))|_{B'}$. Soient $L$ un
corps quadratique sur $K$ et $\sigma $ l'involution qui fixe $K$. Soient $b\in
B(L)$, et  $P_1,\ldots,P_n \in X_b(L)$.  On
dispose d'un point $K$-rationnel $y_{(P_1,\ldots,P_n)}$ de
$Y_n$:
$$y_{(P_1,\ldots,P_n)} = \{(P_1,\ldots,P_n),(P_1,\ldots,P_n)^\sigma\}.$$
 {\bf Terminologie:} Une $K$-vari\'et\'e $X$ est dite {\em de Lang} si
$X_{\overline{K}}$ domine une vari\'et\'e de type g\'en\'eral.
\begin{prp} Suppsons que les fibres de $X\rar B$ soient de genre
$>1$, et soit $B'\subset (B)^2/S_2$ une sous-vari\'et\'e. Il existe un
entier $n>0$ tel que $Y_{n,B'}$ soit une vari\'et\'e de Lang.
\end{prp}
{\bf D\'emonstration:} Voir \cite{chm}, chapitre 4-5: l'argument de
\cite{chm} 5.1 et 5.2 montre qu'il suffit de consid\'erer le cas suivant: soit
$X\rar B$
une famille de variation maximale ayant
une action d'un groupe fini $G$, on doit montrer que pour $n$ assez grand, la
vari\'et\'e $Sym^2(X^n_B)/ G$ est de type g\'en\'eral; ce qui est
d\'emontr\'e  dans \cite{chm} 4.3.

[B. Hassett \cite{hassett} a d\'emontr\'e un r\'esultat plus
fort, en utilisant la th\'eorie  des surfaces stables de Koll\'ar
et Shepherd-Barron: une puissance relative assez grande d'une famille
de surfaces de type g\'en\'eral est une vari\'et\'e de Lang.]

Vu la proposition 1, la conjecture de Lang dit que l'ensemble des points
rationnels de $Y_{n,B'}$
n'est pas dense. Donc nous devons contr\^oler les points d'une
sous-vari\'et\'e $D\subset Y_{n,B'}$, qui poss\`ede un morphisme $D \rar
Y_{n-1}$ de dimension relative $\leq 1$. La proposition suivante nous y aidera:
\begin{prp} Soit $f:D\rar Z$ une famille plate de dimension relative pure 1.
Supposons que toute composante irr\'eductible d'une fibre g\'eom\'etrique est
de genre
g\'eom\'etrique $\geq 2$. Alors, pour $n$ assez grand, toute composante
irr\'eductible de  $D^n_Z$ est une vari\'et\'e de Lang.
\end{prp}
{\bf D\'emonstration:} Nous allons nous ramener au cas o\`u les fibres
sont irr\-\'e\-duc\-tibles, tra\^{\i}t\'e dans \cite{chm}. On peut supposer
$K=\overline{K}$, $D$
irr\'eductible et normale, et on peut remplacer $Z$ par un ouvert, de sorte
qu'on peut supposer  $f$  lisse.  Soit $D\stackrel{f'}{\rar}
Z'\rar Z$ la factorisation de Stein de $f$; les fibres de $f'$ sont
conn\`exes. Si $z\in Z$ est un point g\'en\'eral,
 posons $D_z= \cup_{1\leq i\leq c} C_i$, o\`u les $C_i$ sont les
composantes conn\`exes de $D_z$. Soit $M$ une composante conn\`exe de
 $D^n_Z$ pour $n$ quelconque. Au-dessus de $z$, la fibre de $M$ est form\'ee de
$n$-tuples de points $(P_j)_{j=1}^n$. Pour tout $i$ il existe un ensemble
d'indices $J_i$ tel que pour $j\in J_i$ on a $P_j\in C_i$. Choisissons
 $J_i$ de cardinal $\geq n/c$. Si $n$ est assez grand, alors $\# J_i=n_i$ l'est
aussi.
Nous disposons d'un morphisme $M\rar D_Z^{J_i}\cong D_Z^{n_i}$. L'image
de $M$ dans $D_Z^{J_i}$ est la sous-vari\'et\'e $D_{Z'}^{J_i}$ des $n_i$-tuples
de points qui sont contenus dans la m\^eme composante d'une fibre de $f$.
D'apr\`es \cite{chm}, si $n_i$ est
assez grand $D_{Z'}^{J_i}$ est une vari\'et\'e de Lang, donc $M$ l'est aussi.

Nous aurons besoin d'un resultat semblable a proposition 2, o\`u $D$ n'est pas
de dimension relative pure 1.  Commen\c{c}ons par un lemme facile.
\begin{lem} Soit $D \rar Z$  un morphisme g\'en\'eriquement fini, et soit
$\Delta_n$ la ``grande diagonale''
de $D^n_Z$ (des $n$-tuples des points tels que au moins deux
entre eux sont \'egaux). Alors il existe un entier $n$ tel que $(D^n_Z
\setminus
\Delta_n) \rar Z$ n'est pas dominant.
 \end{lem}
{\bf D\'emonstration:}  Si le degr\'e g\'en\'erique de $D\rar Z$ est au plus
$n-1$,
alors chaque $n$-tuple des points dans une fibre g\'en\'erale contiens au moins
deux points \'egaux.
\begin{prp} Soit $D \rar Z$ de dimension relative g\'en\'erique 1, tel que
toute composante irr\'eductible de dimension 1
d'une fibre g\'en\'erale est de genre  g\'eo\-m\'e\-t\-rique au moins 2, et
soit $\Delta_n$ la grande diagonale de
$D^n_Z$. Alors pour $n$ assez grand, toute composante
irr\'eductible de  $D^n_Z\setminus \Delta_n$ qui domine $Z$   est une
vari\'et\'e de Lang.
 \end{prp}
{\bf D\'emonstration:} R\'esulte de la proposition 2 et du lemme 1.

{\bf D\'emonstration du th\'eor\`eme}.
La base $B$ de la famille des courbes que nous utiliserons est le sch\'ema de
Hilbert
qui param\'etrise les courbes lisses projectives conn\`exes de genre $g$
plong\'ees
par un syst\`eme
lin\'eaire 3-canonique dans un espace projectif. La famille $X$ est la courbe
universelle au-dessus de $B$.
Chaque courbe  lisse projective conn\`exe de genre $g$ d\'efinie sur un corps
de nombres $L$ est isomorphe \`a une fibre au-dessus d'un point $b\in B(L)$.

{\bf Notation:} Soit $m>n$. Nous disons que $ y_{(P_1,\ldots,P_{n})}\in Y_n(K)
$ est {\em $m$-prolongeable } s'il exist un corps $L, \quad [L:K]=2$, et une
``prolongation'' $y_{(P_1,\ldots, P_{m})}\in Y_m(K)$ tel que pour $1\leq i\leq
m$ on ait
$P_i\in X(L)\setminus X(K)$ et pour $1\leq i \neq j\leq n$ on ait $P_i\neq
P_j$.
 Nous notons $ E_n^m $ l'ensemble des points $m$-prolongeable, et $F_n^m =
\overline{E_n^m}$ l'adh\'erence de Zariski de $E_n^m$.

Nous devons d\'emontrer que $F_n^m = \emptyset$ pour  $m$ assez grand.

\begin{lem} Pour tout $n$ il existe $m(n)$ tel que $F_n^{m(n)+k} =
F_n^{m(n)}$ pour  $k$ positif quelconque.
 \end{lem}
 Pour  $n$ fix\'e, les vari\'et\'es $F_n^m$ sont ferm\'ees
d\'ecroissantes dans l'espace noetherien $Y_n$.

{\bf Notation:} $F_n = F_n^{m(n)}$.
\begin{lem} Soit $I\subset \{1,2,...,n'\}$ un $n$-tuple. La
projection $\pi_I: F_{n'} \rar F_n$ est surjective.
 \end{lem}
{\bf D\'emonstration:} $\pi_I (E^m_{n'}) = E^m_n$.

Nous disposons d'un morphisme naturel $\pi_{n,k}:F_{n+k}\rar
(F_{n+1})_{F_n}^k$.
Notons que $\pi_{n,k}(E^{m(n+k)}_{n+k})\subset (F_{n+1})_{F_n}^k\setminus
\Delta_k$, par d\'efinition. Alors le lemme 1, avec $D=F_{n+1}$ et $Z=F_n$
entra\^{\i}ne le lemme suivant:
\begin{lem} Soit  $y\in F_{n}$. La dimension de la fibre de
$F_{n+1}$  au-dessus de $y$ est $\geq 1$.
 \end{lem}

{\bf Lemme 4bis:} {\em Soit $y\in F_{n}$. La dimension de la fibre de
$F_{n+k}$  au-dessus de $y$ est $\geq k$.}

{\bf D\'emonstration:} Induction sur $k$ et Lemme 4.
\begin{lem} Supposons que la conjecture de Lang soit vraie. Soit $y\in
F_{n}$. Alors la dimension de la fibre
de $F_{n+1}$  au dessus de $y$ est  $2$.
 \end{lem}
{\bf D\'emonstration:} Supposons qu'il existe une composante irr\'eductible de
$F_{n+1}$ dont les
fibres g\'en\'erales sont de dimension $<2$. Le lemme 4 dit que les fibres sont
de dimension 1, donc la dimension relative de $(F_{n+1})_{F_n}^k$ est $k$. De
plus, le lemme 4bis entra\^{\i}ne que la dimension relative de $F_{n+k}$ est
$k$
aussi. Donc il existe au moins une composante irr\'eductile $H_k$ de $F_{n+k}$
qui domine une
composante de  $(F_{n+1})_{F_n}^k$. Notons que, puisque $(P_1,\ldots, P_{n})$
est d\'efinie sur $L$
mais  pas sur $K$,  la fibre g\'eom\'etrique au-dessus de $y_{(P_1,\ldots,
P_{n})}\in Y_n$ dans  $Y_{n+1}$ est un produit de deux courbes
de genre $g$. Donc chaque courbe contenue dans cette fibre
est de genre g\'eom\'etrique au moins $g$. La proposition 3, avec $D=F_{n+1}$
et $Z=F_n$, enta\^{\i}ne que pour $k$ assez grand, la composante $H_k$ de
$F_{n+k}$ est une vari\'et\'e de Lang.
 La conjecture de Lang entra\^{\i}ne que l'ensemble des points rationnels de
$H_k$ n'est pas dense dans $H_k$ pour $k$ assez grand, mais $F_{n+k}$ est
definie comme l'adh\'erence de Zariski de son ensemble de points rationnels.

{\bf Lemme 5bis:} {\em  Supposons que la conjecture de Lang soit vraie. Soit
$y\in F_{n}$. Alors la dimension  de la fibre de $F_{n+k}$  au desus de $y$ est
$2k$.}
\begin{lem}  Supposons que la conjecture de Lang soit vraie. Il existe
une sous-vari\'et\'e
$B'\subset B^2/S_2$ tel que pour $k$ positif, si on ecrit $I = \{
n+1,\ldots,n+k\}$ on a $\pi_I(F_{n+k})=Y_{k,B'}$.
 \end{lem}
{\bf D\'emonstration:} Appliquons le lemme 5bis, et prenons $B'$ \'egal \`a
l'image de
$F_{n+k}$ dans $B^2/S_2$.
\begin{lem} On a $F_n = Y_{n,B'}$. \end{lem}
{\bf D\'emonstration:} Lemmes 3 et 6 (ou bien: prenons n=0 dans le lemme 6).

{\bf Conclusion:} La proposition 1 dit que pour $k$ assez grand,
$Y_{k,B'}=F_k$ est une vari\'et\'e de Lang. Par d\'efinition, l'ensemble des
points rationnels de $F_k$ est dense. Par la conjecture de Lang $F_k$ est vide,
et vu le lemme 3, $F_n$ l'est aussi.

{\bf Remarque:} Pour d\'emontrer le cas des corps de degr\'e $d=3$, on remplace
$Y_n= Sym^2(X^n_B)$ par $Sym^3(X^n_B)$,  et on montre comme ci-dessus que la
dimension des fibres de $F_{n+1}\rar F_n$ est au moins 2. On doit
montrer que la dimension est 3. Pour achever ce but, on peut appliquer le
th\'eor\`eme de Hassett \cite{hassett} sur les familles
de surfaces de type g\'en\'eral, et d\'emontrer que si la dimension
des fibres est 2 et si $k$ est assez grand, $F_{n+k} $ est une vari\'et\'e de
Lang.

\small
 \def\refname{{\small \sc R\'ef\'erences Bibliographiques}}

Department of Mathematics, Boston University\\
111 Cummington Street, Boston, MA 02215 \\
abrmovic@math.bu.edu
\end{document}